\documentclass[12pt,a4paper, twoside]{article}

\usepackage{amsmath,amssymb,amsthm,amsfonts,mathrsfs,amscd,environ}
\usepackage{dsfont}
\usepackage{latexsym,enumerate,color,geometry,extarrows}
\usepackage{xcolor}
\usepackage{verbatim,fancyhdr,enumitem}

 \NewEnviron{ews}{%
\begin{equation}\begin{split}
  \BODY
\end{split}\end{equation}
}

\NewEnviron{ews*}{%
\begin{equation*}\begin{split}
  \BODY
\end{split}\end{equation*}
}

 \newtheorem{thm}{Theorem}[section]
 \newtheorem{lem}[thm]{Lemma}
 
 \newtheorem{exa}[thm]{Example}
 
 \newtheorem{defn}{Definition}[section]
 \newtheorem{rem}{Remark}[section]
 
 \numberwithin{equation}{section}

\def\C{{\mathscr{C}}}
\def\trace{{\rm trace}}

\def\dif{{\mathord{{\rm d}}}}
\def\no{\nonumber}

\def\mR{{\mathbb R}}
\def\mE{{\mathbb E}}

\def\sF{{\mathscr F}}

\def\bd{\begin{defn}}
\def\ed{\end{defn}}
\def\bl{\begin{lem}}
\def\el{\end{lem}}
\def\bt{\begin{thm}}
\def\et{\end{thm}}
\def\br{\begin{rem}}
\def\er{\end{rem}}

\allowdisplaybreaks

\title{{\bf Stability of   hybrid   pantograph stochastic functional differential  equations }
}

\author{
{\bf Hao Wu$^{a)}$, Junhao Hu$^{a)}$, Chenggui Yuan$^{b)}$}\\
\footnotesize{$^{a)}$School of Mathematics and Statistics,
South-Central University For Nationalities}\\
\footnotesize{ Wuhan, Hubei 430000, P.R.China}\\
\footnotesize{Email: wuhaomoonsky@163.com},
\footnotesize{ junhaohu74@163.com}\\
\footnotesize{$^{b)}$Mathematic department, Swansea University, Bay campus, SA1 8EN, UK}\\
\footnotesize{Email: C.Yuan@Swansea.ac.uk}\\
}

\begin{document}

\maketitle

\begin{abstract}
In this paper, we study a new type of stochastic functional differential equations which is called hybrid pantograph stochastic functional differential  equations.   We investigate  several  moment properties and sample properties of the solutions to  the equations by using the method of multiple Lyapunov functions, such as the moment exponential stability,  almost sure exponential stability  and almost sure polynomial stability, etc.
\end{abstract}\noindent

AMS Subject Classification (2020): \quad 60H10; \quad 34K20.
\noindent

Keywords: Moment properties;  Sample properties;   Markovian switching; Pantograph stochastic functional differential  equations

\section{Introduction}
 Stochastic differential equations (SDEs)   are widely used to model stochastic systems  in  different branches of science and industry.  Stability  and boundedness  of the solution   are  the most popular topics in the area of stochastic systems and control.  We refer the reader to \cite{BHB, FP1, CFAF} and references therein.   Dynamic systems may not only depend on present states but also the past states.  Stochastic delay differential equations (SDDEs) and pantograph stochastic delay differential equations (PSDDEs) are often used to model these systems,  whose systems  depend on the past state $x(t-\tau)$ and $x(\theta t)$   respectively. The form of these equations are as follows:
  \begin{align*}
 \dif x(t)= f( x(t), x(t-\tau), t)\dif t +g( x(t),x(t-\tau),t)\dif  B(t),
 \end{align*}
 and
\begin{align*}
 \dif x(t)= f( x(t), x(\theta t), t)\dif t +g( x(t),x(\theta t),t)\dif  B(t),
 \end{align*}
   where $\tau ,\theta $ are two constants satisfying $\tau>0, 0<\theta <1.$ We here only mention \cite{Pengsg, PSG, KKMM, RX1, LR111}, to name a few. However, there are many practical systems whose future state depends on the states over the whole time interval $[t-\tau, t]$ rather than at times $t-\tau$ and $t.$   Stochastic functional differential  equations (SFDEs) have therefore been developed to describe such systems. Generally speaking, SFDEs have the form:
\begin{align*}
 &\dif x(t)= f( x_t, t)\dif t +g(  x_t, t)\dif  B(t), t\in [t_{0}, \infty), \quad x_{t_{0}}=\xi,
 \end{align*}
where $x_{t}=\{x( t-\theta), 0\leq \theta \leq \tau)\},$ $\tau>0$ is a constant. As is well known, many scholars realized that numerous system in our real world  may  experience abrupt changes in their structure and parameters caused by  phenomena  such as component failures or repairs, changing subsystem interconnections and abrupt environmental disturbances. Hybrid systems driven by continuous-time Markov chains have been used to cope with such situation. Markov chains play the role of stabilizing factor in the stability of hybrid systems. That means, when some subsystems are unstable, but others are stable, then the overall system could be stable because of switching between the subsystems.  Since then, the literature on the topic of stability for  stochastic differential equations with Markovian switching (SDEswMS) bloomed, both in the direction of obtaining qualitative and quantitative results for the generalized emerging equations and on developing applications which aim to  population ecology, network, heat exchanges, etc.    For example, \cite{BBG} studied the stability of semi-linear SDEswMS, and \cite{ES}
 investigated the  following general   nonlinear SDEswMS:
\begin{align*}
 \dif x(t)= f( x(t), t, r(t))\dif t +g( x(t),t,r(t))\dif  B(t),
 \end{align*}
where $r(t)$ is a Markov chain taking values in $S=\{1,2,\cdots, N \}.$
Moreover,  \cite{MMXX4}  applied  SDEswMS to solving a control problems,   \cite{MMXX5} investigated some complex-valued coupled oscillators,   \cite{MMX1} studied the stability of regime-switching jump diffusion systems,
\cite{MMXX6} analyzed  asymptotic stability in distribution for such type of equations.
Later,   the study of stochastic functional differential equations with Markovian switching (SFDEswMS) (including stochastic delay differential equations with Markovian switching) and  PSDDEs with Markovian switching  have been  also developed rapidly.
 Many scholars have enthusiastically studied the stability of such equations and given some applications. For example,   \cite{SFM}  investigated the exponential stability of highly nonlinear neutral pantograph stochastic differential equations,   \cite{SR}  built Razumikhin-type theorems on neutral SFDEs,  \cite{BJ} studied  stability of neutral SFDEswMS driven by G-Brownian motion, \cite{HJsS}  analyzed  asymptotic stability and boundedness of SFDESwMS.  More related work can be seen in  \cite{OPZ, glhj, LLM,  LR10, AAA1,  LR111, MMXX8, MMX21}.

To the best of our knowledge, so far  there is  little study  on  hybrid   pantograph stochastic  functional differential  equations (HPSFDEs), in which the $\theta$ changes in interval (0,1], while  the $\theta$  is a constant  in pantograph stochastic delay differential  equations. Inspired by the works of the above articles, we aim in this paper to study several  moment properties and sample properties of the solutions  such as the moment exponential stability,  almost surely exponential stability  and almost sure polynomial stability, etc. for HPSFDEs.

We close this part by giving our organization  in this article. In Section 2, we introduce some necessary notations. In Section 3, we give our main results on  the moment properties and sample properties of  analytical solution. Several  examples are also given  to illustrate the theory.

 \section{Preliminaries}
\subsection{Notations}Throughout this paper,  Let  $(\Omega, \sF, \{\sF_{t}\}_{t\geq 0}, P)$ be  a complete probability space satisfying the usual conditions(i.e., it is increasing and right continuous with $\sF_{0}$ contains all $P$-null sets) taking along
 a standard  $d$-Brownian motion process  $B(t).$  For $x, y \in \mR^{n},$ we use $|x |$  to denote the Euclidean norm of
$x,$ and use  $\langle x, y\rangle$ or $x^{T}y$ to denote the Euclidean inner product. If  $A$ is a matrix, $A^{T}$  is the transpose  of $A$  and   $|A |$ represents  $\sqrt{\mathrm{Tr} (AA^{T}).}$   Let $\lfloor a \rfloor$ be the integer parts of $a.$  Moreover, for $0<\underline{\theta} <1,$ denote by $\mathscr{C} :=\mathscr{C}([\underline{\theta}, 1]; \mR^{n})$ the family of all continuous $\mR^{n}-$valued functions $\varphi$  defined on $[\underline{\theta}, 1]$ with the norm $\| \varphi\|=\sup_{\underline{\theta} \leq t\leq 1}|\varphi(t)|.$ Let  $t_0>0$ and $h: [t_0, \infty)\to \mR^{n}$  be a continuous function,   for $t \ge t_0$ denote $h_t(\theta)=h(\theta t),  \underline{\theta} \le \theta \le 1.$  One can see that $h_t(\cdot)\in \mathscr{C}.$
Let $r(t)$ be a   continuous-time Markov chain taking values in $S=\{1,2,\cdots, N\}$ with the generator $\Gamma=(\gamma_{ij})_{N\times N}$  such that
\begin{align*}
 P\{r(t+\delta)=j|r(t)=i\}=
 \begin{cases}
 \gamma_{ij}\delta + o(\delta), & i\neq j, \\
 1+\gamma_{ii}\delta + o(\delta), & i=j,
 \end{cases}
\end{align*}
where $\delta > 0.$ Here $\gamma_{ij} $ is the transition rates from $i$ to $j$ and $\gamma_{ij}\ge 0$ if $i\neq j$ while $\gamma_{ii}=-\sum_{j\neq i}\gamma_{ij}.$ It is  well know that almost every sample path of $r(t)$ is a right-continuous step functions with finite number of sample jumps in any finite subinterval of $\mR_{+}=[0, \infty).$  Assume that Markov chain $r(t)$ is independent of Brownian motion.

Denote by $ C^{1,2}([t_0,+\infty)\times R^{n} \times S; [0,+\infty))$ the
  family of all continuous nonnegative functions $V(t,x, i)$ defined on $[t_0,+\infty)\times R^{n} \times S $,
  such that for each $i\in S$,
   they are continuously once  differentiable in $t$ and twice in $x$.

\section{Main Results}
 Consider the following  HPSFDE:
\begin{equation}\label{eq1}
\begin{cases}
 &\dif x(t)= f( x_t, t,r(t))\dif t +g(  x_t, t, r(t))\dif  B(t), t\in [t_{0}, \infty),\\
 & x(t)=\xi(t), t\in [\underline{\theta} t_{0}, t_{0}],
 \end{cases}
 \end{equation}
where $x_{t}=\{x(\theta t), \underline{\theta}\leq \theta \leq 1)\}$ and $0< \underline{\theta}<1$ is a constant. We would like to point out that $x_{t}\in \C$ is a segment process and $x_t(\theta)=x(\theta t)$ while $x(t)\in R^n$ is a point.

Given $V\in C^{1,2}( R^{n} \times [t_0,+\infty)\times S; [0,+\infty)), \varphi\in \C, i\in S$,
    we define an operator  $LV: \mathscr{C}\times [t_{0}, \infty)\times S  \rightarrow \mR $ by
\begin{align*}
LV(\varphi,t,i)&=V_{t}(\varphi(1),t,i)+ V_{x}(\varphi(1),t,i)f(\varphi,t,i)
&+\frac{1}{2}\trace(g^{T}(\varphi,t,i)V_{xx}(\varphi(1),t,i)g(\varphi,t,i))\\
&+\sum^{N}_{l=1}\gamma_{il}V(\varphi(1),t,l),
 \end{align*}
where
\begin{align*}
V_t(x, t, i)=\left(\frac{\partial V(x, t, i)}{\partial t}\right), \, V_x(x, t, i)=\left(\frac{\partial V(x, t, i)}{\partial x_1}, \ldots, \frac{\partial V(x, t, i)}{\partial x_n}\right)
\end{align*}
and
\begin{align*}
 V_{xx}(x, t, i)=\left(\frac{\partial^2 V(x, t, i)}{\partial x_k \partial x_l}\right)_{kl}.
\end{align*}
 We have the following  the corresponding It\^{o}'s formula for hybrid system \eqref{eq1}:
\\
\begin{align*}
&V(x(t),t,r(t))
=V(x(0),0,r(0))+ \int^{t}_{t_0}LV(x_s, s,r(s))\dif s\\
 &+ \int^{t}_{t_0}V_{x}(x(s),s,r(s))g(x_{s},s,r(s))\dif B(s).
 \end{align*}

 The following assumptions are needed.
\begin{itemize}
\item[(H1)] For any $\varphi, \varphi' \in C([\underline{\theta}, 1]; \mR^{n})$ satisfying $\|\varphi\|\vee \|\varphi'\|\leq R,$  there exists a positive constant $C_{R}$ such that
    \begin{align*}
    |f( \varphi, t,i) - f( \varphi^{\prime},t,i)|\vee|g( \varphi, t,i) - g( \varphi^{\prime},t,i)|
    \leq C_{R}\|\varphi- \varphi^{\prime}\|.
    \end{align*}
\item[(H2)]  There exist functions $V\in C^{2,1}(\mR^{n}\times [t_{0}, \infty)\times S; \mR_{+}),$ $U_{0}, U_{k}\in C^{2,1}(\mR^{n}\times [t_{0}, \infty); \mR_{+}),$ and probability measures $\nu_{k}$ on $[\underline{\theta}, 1]$,     and  non-negative constants $a_{0}, a_{k}, b_{kl} ,$ $ k=1,2,\cdots, M,$ $ l=1,2,\cdots, l_{k} $  such that
  \begin{align}
 \lim_{|x|\rightarrow \infty}\inf_{t_{0}\leq t< \infty}U_{0}(x,t)=\infty,
 \end{align}
\begin{align}
U_{0}(x,t)\leq V(x,t,i)\leq U_{1}(x,t), \forall (x,t,i)\in \mR^{n}\times \mR^{+}\times S,
 \end{align}
\begin{align}\label{L0}
LV(\varphi,t,i)&\leq a_{0}+\sum^{M}_{k=1}\bigg[\bigg.-a_{k}U_{k}(\varphi(1), t)\nonumber\\
&+ \sum^{l_{k}}_{l=1}b_{kl}\int^{1}_{\underline{\theta}} e^{-\int^{t}_{0}\lambda(\theta, u)\dif u}U_{k}^{\alpha_{kl}}(\varphi(1), t) U_{k}^{1-\alpha_{kl}}(\varphi(\theta), \theta t)\dif \nu_{k}(\theta) \bigg]\bigg.,
 \end{align}
where  function $\lambda(\cdot, \cdot):[\underline{\theta}, 1]\times \mR_{+}\rightarrow \mR_{+}$ satisfying $\inf_{0 \leq s < \infty}\lambda (\theta, s)\geq \beta (1-\theta)$  and $\alpha_{kl}, \beta$ are   constants satisfying  $0\leq\alpha_{kl}\leq 1, 0<\beta<a_{1}$.
\end{itemize}

\subsection{Existence and Uniqueness}

In the same way as in \cite{LR1211},  we can show that \eqref{eq1}   has  a unique local solution $x(t), t\in [t_0, \sigma_{\infty})$  under $(\mathrm{H1}),$  where $\sigma_{\infty}$ is the explosion time.  The following condition (H2') will guarantee  a global solution to \eqref{eq1}, that is
\begin{itemize}
\item[(H2')]  Assume that (H2) holds, but \eqref{L0} is replaced by
\begin{align}
LV(\varphi,t,i)&\leq a_{0}+\sum^{M}_{k=1}\bigg[\bigg.-a_{k}U_{k}(\varphi(1), t)\nonumber\\
&+ \sum^{l_{k}}_{l=1}b_{kl}\int^{1}_{\underline{\theta}} U_{k}^{\alpha_{kl}}(\varphi(1), t) U_{k}^{1-\alpha_{kl}}(\varphi(\theta), \theta t)\dif \nu_{k}(\theta) \bigg]\bigg.,
 \end{align}
where  $0\leq \alpha_{kl}\leq 1$ are   constants .
\end{itemize}

\bt
Assume that (\rm{H1}) and (\rm{H2'}) hold.  If
\begin{align}\label{eu1}
-a_{k}+ \sum^{l_{k}}_{l=1}b_{kl}\alpha_{kl}
+\sum^{l_{k}}_{l=1}b_{kl}\frac{1}{\underline{\theta}}(1-\alpha_{kl})\le 0, \, k=1,2,\cdots, M,
\end{align}
  then the equation \eqref{eq1} has a unique global solution.
\et
\begin{proof}
Let $x(t), t\in [t_0, \sigma_{\infty})$ be the unique local solution and $\sigma_{a}=\inf\{t\geq t_0: |x(t)|\geq a\}.$
Using It\^{o}'s formula and taking the expectation, we have
\begin{align}
&\mE[V(x(t\wedge \sigma_a), t\wedge \sigma_{a}, r(t\wedge \sigma_a))]=\mE[V(x(t_0), t_0, r(t_0))] + \mE\int^{t\wedge  \sigma_a}_{t_{0}}LV(x_{s},s,r(s))\dif s  \no\\
& \leq \mE[V(x(t_0), t_0, r(t_0))] +  \mE\int^{t\wedge  \sigma_a}_{t_{0}}\bigg\{\bigg. a_{0}+\sum^{M}_{k=1}\bigg[\bigg.-a_{k}U_{k}(x(s), s)\no\\
& \quad \quad \quad \quad \quad \quad\quad \quad \quad+ \sum^{l_{k}}_{l=1}b_{kl}\int^{1}_{\underline{\theta}} U_{k}^{\alpha_{kl}}(x(s), s) U_{k}^{1-\alpha_{kl}}(x(\theta s), \theta s)\dif \nu_{k}(\theta) \bigg]\bigg.\bigg\}\bigg.\dif s\no\\
&\leq \mE[V(x(t_0), t_0, r(t_0))] +  \mE\int^{t\wedge  \sigma_a}_{t_{0}}\bigg\{\bigg. a_{0}+\sum^{M}_{k=1}\bigg[\bigg.-a_{k}U_{k}(x(s), s)\no\\
& \quad \quad +  \sum^{l_{k}}_{l=1}b_{kl}\alpha_{kl}U_{k}(x(s), s)+ \sum^{l_{k}}_{l=1}b_{kl}(1-\alpha_{kl})\int^{1}_{\underline{\theta}} U_{k}(x(\theta s), \theta s)\dif \nu_{k}(\theta) \bigg]\bigg.\bigg\}\bigg.\dif s.
 \end{align}
 Noting that
 \begin{equation*}
 \begin{split}
 &\int^{t\wedge  \sigma_a}_{t_{0}}\int^{1}_{\underline{\theta}}U_k(x(\theta s), \theta s)\dif \nu_k(\theta)\dif s =\int^{1}_{\underline{\theta}}\int^{t\wedge  \sigma_a}_{t_{0}}U_k(x(\theta s), \theta s)\dif s\dif \nu_k(\theta)\\
 &\le \frac{1}{{\underline{\theta}}}\int^{1}_{\underline{\theta}}\int^{\theta(t\wedge  \sigma_a)}_{\theta t_{0}}U_k(x( s),  s)\dif s\dif \nu_k(\theta)\\
 &\le  \frac{1}{{\underline{\theta}}}\int^{1}_{\underline{\theta}}\int^{t\wedge  \sigma_a}_{\underline{\theta} t_{0}}U_k(x( s),  s)\dif s\dif \nu_k(\theta)\\
 &\le \frac{1}{{\underline{\theta}}}\int^{1}_{\underline{\theta}}\int^{t\wedge  \sigma_a}_{t_{0}}U_k(x( s),  s)\dif s\dif \nu_k(\theta)+  \frac{1}{{\underline{\theta}}}\int^{1}_{\underline{\theta}}\int^{t_0}_{\underline{\theta} t_{0}}U_k(x( s),  s)\dif s\dif \nu_k(\theta)\\
 &\le \frac{1}{{\underline{\theta}}}\int^{t\wedge  \sigma_a}_{t_{0}}U_k(x( s),  s)\dif s +  \frac{1}{{\underline{\theta}}}\int^{t_0}_{\underline{\theta} t_{0}}U_k(x( s),  s)\dif s,
 \end{split}
 \end{equation*}
  and \eqref{eu1}, one can see that
 \begin{align}\label{3.8}
&\mE[V(x(t\wedge \sigma_a), t\wedge \sigma_{a}, r(t\wedge \sigma_a))]\no\\
&\leq \mE[V(x(t_0), t_0, r(t_0))] +  \mE\int^{t\wedge  \sigma_a}_{t_{0}}\bigg\{\bigg. a_{0}+\sum^{M}_{k=1}\bigg[\bigg.-a_{k}U_{k}(x(s), s)\no\\
&  \quad \quad \quad\quad \quad +  \sum^{l_{k}}_{l=1}b_{kl}\alpha_{kl} U_{k}(x(s), s)+ \sum^{l_{k}}_{l=1}b_{kl}\frac{1}{\underline{\theta}}(1-\alpha_{kl})U_{k}(x( s),  s) \bigg]\bigg.\bigg\}\bigg.\dif s\no\\
& + \sum^{M}_{k=1}\sum^{l_{k}}_{l=1}b_{kl}\frac{1}{\underline{\theta}}(1-\alpha_{kl})\int^{t_{0}}_{\underline{\theta}t_{0}}U_{k}(x( s),  s)\dif s                          \no\\
&\leq \mE[V(x(t_0), t_0, r(t_0))] \no\\
&+  \mE\int^{t\wedge  \sigma_a}_{t_{0}}\bigg\{\bigg. a_{0}+\sum^{M}_{k=1}\bigg(\bigg.-a_{k}+ \sum^{l_{k}}_{l=1}b_{kl}\alpha_{kl}
+\sum_{k=1}^{l_{k}} b_{kl}\frac{1}{\underline{\theta}}(1-\alpha_{kl})\bigg)\bigg.U_{k}(x(s), s)\bigg\}\bigg.\dif s\no\\
& + \sum^{M}_{k=1}\sum^{l_{k}}_{l=1}b_{kl}\frac{1}{\underline{\theta}}(1-\alpha_{kl})\int^{t_{0}}_{\underline{\theta}t_{0}}\mE U_{k}(\xi( s),  s)\dif s \no\\
 &\leq c_{0} +a_{0}t,
 \end{align}
where $c_{0}=\mE[V(x(t_0), t_0, r(t_0))]+ \sum^{M}_{k=1} \sum^{l_{k}}_{l=1}b_{kl}\frac{1}{\underline{\theta}}(1-\alpha_{kl})\int^{t_{0}}_{\underline{\theta}t_{0}}\mE U_{k}(\xi( s),  s)\dif s .$

Setting $\mu_{a}= \inf_{|x|\geq a, t_{0}\leq t< \infty}U_{0}(x, t),$ we then have
\begin{align*}
\mE[U_{0}(x(t\wedge  \sigma_a), t\wedge  \sigma_a)]\geq \mE[U_{0}(x( \sigma_a),  \sigma_a)1_{\sigma_{a}\leq t}]\geq \mu_{a}P(\sigma_{a}\leq t).
\end{align*}
This immediately implies
\begin{align*}
P(\sigma_{\infty}\leq t)=\lim_{a\rightarrow \infty} P(\sigma_{a}\leq t)\leq \lim_{a\rightarrow \infty} \frac{\mE[U_{0}(x(t\wedge  \sigma_a), t\wedge \sigma_{a}]}{\mu_{a}}=\lim_{a\rightarrow \infty}\frac{c_{0} +a_{0}t}{\mu_{a}}=0.
\end{align*}
Therefore,  $\sigma_{\infty}=\infty,\, a.s.,$ and  there exists unique global  solution $x(t)$ on $[t_{0}, \infty).$
 \end{proof}

 \subsection{Exponential Stability}

In this subsection, we will investigate several  moment properties and sample properties of the solutions to  the equations  such as the moment exponential stability,  almost sure exponential stability, etc. Before studying the stability of the solution to E.q.\eqref{eq1}, we present a semi-martingale convergence theorem which can be found in \cite{LSLS}.
\bl
Let $A_{1}(t), A_{2}(t)$ be two continuous adapted increasing processes on $t\geq 0$ with $A_{1}(0)=0, A_{2}(0)=0,$ a.s., M(t) a real-valued continuous local martingale with $M(0)=0,$ a.s., $\xi$ a nonnegative $\sF_{0}-$measurable random variable such that $\mE[\xi]< \infty.$ Set $X(t)=\xi+A_{1}(t)-A_{2}(t)+M(t), t\geq 0.$  If $X(t)$ is nonnegative , then we have the following results:
$$\{\lim_{t\rightarrow \infty}A_{1}(t)<\infty\}\subset \{\lim_{t\rightarrow \infty}A_{2}(t)<\infty\}\cap \{\lim_{t\rightarrow \infty}X(t)<\infty\},\,a.s.,$$
where $C\subset D,\, a.s.$  means $P(C\cap D^{c})=0.$ In particular, if $\lim_{t\rightarrow \infty}A_{1}(t)<\infty,\,a.s.,$ then, with probability one,
$$\lim_{t\rightarrow \infty}A_{2}(t)<\infty, \,  \lim_{t\rightarrow \infty}X(t)<\infty, \, -\infty <\lim_{t\rightarrow \infty}M(t)<\infty    ,\,a.s.$$
\el

 \bt\label{T3.3}
Assume that $(\mathrm{\mathrm{H}1})-(\mathrm{H2})$ hold with  $$-a_{k}+ \sum^{l_{k}}_{l=1}b_{kl}\alpha_{kl}
+\sum^{l_{k}}_{l=1}b_{kl}\frac{1}{\underline{\theta}}(1-\alpha_{kl})<0, k=1,2,\cdots, M.$$   We then have the following results:
\begin{itemize}
\item[{\rm (i)}] $\limsup_{t\rightarrow \infty}\mE[U_{0}(x(t_{0},t,\xi,i_{0}), t)]\leq \frac{a_{0}}{\varepsilon},$
where $0<\varepsilon \leq \beta$ is a constant satisfying $$a_{1}-\varepsilon- \sum^{l_{1}}_{l=1}b_{1l}\alpha_{1l}- \sum^{l_{1}}_{l=1}b_{1l}\frac{1}{\underline{\theta}}(1-\alpha_{1l})>0.$$\\
\item[{\rm (ii)}]  \begin{align*}
&\limsup_{t\rightarrow \infty}\frac{1}{t}\int^{t}_{t_{0}}\mE[U_{k}(x(t_{0},s,\xi,i_{0}), s)]\dif s \\
&\leq \frac{a_{0}}{a_{k}-\sum^{l_{k}}_{l=1}b_{kl}e^{-\beta(1-\underline{\theta})t_{0}}\alpha_{kl}
 -\sum^{l_{k}}_{l=1}b_{kl}\frac{1}{\underline{\theta}}e^{-\beta(1-\underline{\theta})t_{0}}(1-\alpha_{kl})}, \quad
  k=1,2,\cdots M.
  \end{align*}
\item[{\rm (iii)}] If $a_{0}=0, $  then the global solution $x(t_{0},t,\xi,i_{0})$ is exponentially stable in moment and almost surely exponential  stable, i.e.
\begin{align}\label{3c}
\limsup_{t\rightarrow \infty}\frac{1}{t}\log (\mE[U_{0}(x(t_{0},t,\xi,i_{0}), t)]) \leq -\varepsilon.
 \end{align}
\begin{align}\label{3c1}
\limsup_{t\rightarrow \infty}\frac{1}{t}\log (U_{0}(x(t_{0},t,\xi,i_{0}), t)) \leq -\varepsilon,\,a.s.,
 \end{align}
 where $\varepsilon$ satisfies the condition in (i).
\end{itemize}
\et
\begin{proof}
(i) Using It\^{o}'s formula  to $e^{\varepsilon t}V(x(t), t, r(t))$, we have
\begin{equation}\label{f1}
\begin{split}
&\mE[e^{\varepsilon (t\wedge \sigma_{a})}V(x(t\wedge \sigma_a), t\wedge \sigma_a, r(t\wedge \sigma_a))]\\
&=\mE[e^{\varepsilon t_{0}}V(x(t_0), t_0, r(t_0))] + \mE\int^{t\wedge  \sigma_a}_{t_{0}}e^{\varepsilon s}(\varepsilon V(x_{s},s,r(s)) +LV(x_{s},s,r(s)))\dif s  \\
& \leq \mE[e^{\varepsilon t_{0}}V(x(t_0), t_0, r(t_0))] \\
& \quad \quad +  \mE\int^{t\wedge  \sigma_a}_{t_{0}}e^{\varepsilon s}\bigg\{\bigg. a_{0}-(a_{1}-\varepsilon)U_{1}(x(s), s)+\bigg[\bigg.-\sum^{M}_{k=2}a_{k}U_{k}(x(s), s)\\
&  \quad \quad+\sum^{M}_{k=1} \sum^{l_{k}}_{l=1}b_{kl}\int^{1}_{\underline{\theta}} e^{-\int^{s}_{0}\lambda(\theta, u)\dif u}U_{k}^{\alpha_{kl}}(x(s), s) U_{k}^{1-\alpha_{kl}}(x(\theta s), \theta s)\dif \nu_{k}(\theta) \bigg]\bigg.\bigg\}\bigg.\dif s.
\end{split}
\end{equation}
Now, we compute
\begin{align}\label{3.13}
 &b_{kl}\int^{t\wedge \sigma_{a}}_{t_{0}}\int^{1}_{\underline{\theta}} e^{\varepsilon s-\int^{s}_{0}\lambda(\theta, u)\dif u}U_{k}^{\alpha_{kl}}(x(s), s) U_{k}^{1-\alpha_{kl}}(x(\theta s), \theta s)\dif \nu_{k}(\theta)\dif s\no\\
 & \leq b_{kl}\alpha_{kl} \int^{t\wedge \sigma_{a}}_{t_{0}}\int^{1}_{\underline{\theta}}e^{\varepsilon s-\int^{s}_{0}\lambda(\theta, u)\dif u}U_{k}(x(s), s) \dif \nu_{k}(\theta) \dif s\no\\
 &+b_{kl}(1-\alpha_{kl})\int^{t\wedge \sigma_{a}}_{t_{0}}\int^{1}_{\underline{\theta}} e^{\varepsilon s-\int^{s}_{0}\lambda(\theta, u)\dif u} U_{k}(x(\theta s), \theta s)\dif \nu_{k}(\theta) \dif s\no\\
 & \leq b_{kl}\alpha_{kl}\int^{1}_{\underline{\theta}}e^{-\beta(1-\theta)t_{0}}\dif \nu_{k}(\theta) \int^{t\wedge \sigma_{a}}_{t_{0}}e^{\varepsilon s}U_{k}(x(s), s)  \dif s\no\\
 &+b_{kl}(1-\alpha_{kl})\int^{t\wedge \sigma_{a}}_{t_{0}}\int^{1}_{\underline{\theta}} e^{\varepsilon s-\int^{s}_{0}\lambda(\theta, u)\dif u} U_{k}(x(\theta s), \theta s)\dif \nu_{k}(\theta) \dif s\no\\
&\leq b_{kl}\alpha_{kl}\int^{t\wedge \sigma_{a}}_{t_{0}}e^{\varepsilon s}U_{k}(x(s), s) \dif s+b_{kl}(1-\alpha_{kl})\int^{1}_{\underline{\theta}}\int^{t\wedge \sigma_{a}}_{t_{0}} e^{\varepsilon s-\int^{s}_{0}\lambda(\theta, u)\dif u} U_{k}(x(\theta s), \theta s)\dif s\dif \nu_{k}(\theta) \no\\
&\leq b_{kl}\alpha_{kl}\int^{t\wedge \sigma_{a}}_{t_{0}}e^{\varepsilon s}U_{k}(x(s), s) \dif s+b_{kl}\frac{1}{\underline{\theta}}(1-\alpha_{kl})\int^{1}_{\underline{\theta}}\int^{t\wedge \sigma_{a}}_{\underline{\theta}t_{0}} e^{\frac{\varepsilon }{\theta}s-\int^{\frac{s}{\theta}}_{0}\lambda(\theta, u)\dif u} U_{k}(x(s),  s)\dif s\dif \nu_{k}(\theta) \no\\
&\leq b_{kl}\alpha_{kl}\int^{t\wedge \sigma_{a}}_{t_{0}}e^{\varepsilon s}U_{k}(x(s), s) \dif s+  b_{kl}\frac{1}{\underline{\theta}}(1-\alpha_{kl})\int^{t\wedge \sigma_{a}}_{t_{0}} e^{\varepsilon s} U_{k}(x( s),  s)\dif s \no\\
& +  b_{kl}\frac{1}{\underline{\theta}}(1-\alpha_{kl})\int^{t_{0}}_{\underline{\theta}t_{0}} e^{\varepsilon s} U_{k}(x( s),  s)\dif s.
\end{align}
This, together with \eqref{f1}, yields that
\begin{equation}\label{f2}
\begin{split}
&\mE[e^{\varepsilon (t\wedge \sigma_{a})}V(x(t\wedge \sigma_a), t\wedge \sigma_a, r(t\wedge \sigma_a))]\\
& \leq \bar{c}_{0}
+  \mE\int^{t\wedge  \sigma_a}_{t_{0}}e^{\varepsilon s}\bigg\{\bigg. a_{0}-\bigg(\bigg.a_{1}-\varepsilon- \sum^{l_{1}}_{l=1}b_{1l}\alpha_{1l}- \sum^{l_{1}}_{l=1}b_{1l}\frac{1}{\underline{\theta}}(1-\alpha_{1l})\bigg)\bigg.U_{1}(x(s), s)\\
&\quad \quad + \sum^{M}_{k=2}\bigg(\bigg.-a_{k}+\sum^{l_{k}}_{l=1}b_{kl}\alpha_{kl}+\sum^{l_{k}}_{l=1}b_{kl}\frac{1}{\underline{\theta}}(1-\alpha_{kl})\bigg)\bigg.U_{k}(x(s), s)\bigg\}\bigg. \dif s \\
& \leq \bar{c}_{0} +  \frac{a_{0}}{\varepsilon}e^{\varepsilon t},
\end{split}
\end{equation}
where $\bar{c}_{0}=\mE[e^{\varepsilon t_{0}}V(x(t_0), t_0, r(t_0))]+   \sum^{M}_{k=1}\sum^{l_{k}}_{l=1}\mE\int^{t_0}_{\underline{\theta}t_{0}}e^{\varepsilon s}\frac{1}{\underline{\theta}}b_{kl}(1-\alpha_{kl})U_{k}(\xi(s), s) \dif s.$
Letting $a \rightarrow \infty,$  it leads to
\begin{align*}
\mE[e^{\varepsilon t}U_{0}(x(t), t, r(t))]
 \leq \bar{c}_{0} +  \frac{a_{0}}{\varepsilon}e^{\varepsilon t}.
\end{align*}
The assertion (i) follows by letting $t \rightarrow \infty$.

(ii) Similar to the proofs of \eqref{3.8} and \eqref{3.13}, we can show that
\begin{equation}
\begin{split}
&\mE[V(x(t\wedge \sigma_a), t\wedge \sigma_{a}, r(t\wedge \sigma_a))]=\mE[V(x(t_0), t_0, r(t_0))] + \mE\int^{t\wedge  \sigma_a}_{t_{0}}LV(x_{s},s,r(s))\dif s  \\
& \leq \mE[V(x(t_0), t_0, r(t_0))] +  \mE\int^{t\wedge  \sigma_a}_{t_{0}}\bigg\{\bigg. a_{0}+\sum^{M}_{k=1}\bigg[\bigg.-a_{k}U_{k}(x(s), s)\\
& \quad \quad \quad+ \sum^{l_{k}}_{l=1}b_{kl}\int^{1}_{\underline{\theta}} e^{-\int^{t}_{0}\lambda(\theta, u)\dif u}U_{k}^{\alpha_{kl}}(x(s), s) U_{k}^{1-\alpha_{kl}}(x(\theta s), \theta s)\dif \nu_{k}(\theta) \bigg]\bigg.\bigg\}\bigg.\dif s\\
&\leq \mE[V(x(t_0), t_0, r(t_0))] \\
&+  \mE\int^{t\wedge  \sigma_a}_{t_{0}}\bigg\{\bigg. a_{0}+\sum^{M}_{k=1}\bigg(\bigg.-a_{k}+ \sum^{l_{k}}_{l=1}b_{kl}e^{-\beta(1-\underline{\theta})t_{0}}\alpha_{kl}\\
& \quad \quad \quad \quad\quad \quad \quad \quad\quad \quad \quad \quad+\sum_{k=1}^{l_{k}} b_{kl}\frac{1}{\underline{\theta}}e^{-\beta(1-\underline{\theta})t_{0}}(1-\alpha_{kl})\bigg)\bigg.U_{k}(x(s), s)\bigg\}\bigg.\dif s\\
& + \sum^{M}_{k=1}\sum^{l_{k}}_{l=1}b_{kl}\frac{1}{\underline{\theta}}e^{-\beta(1-\underline{\theta})t_{0}}(1-\alpha_{kl})\int^{t_{0}}_{\underline{\theta}t_{0}}U_{k}(\xi( s),  s)\dif s \\
 &\leq c_{0} +a_{0}t.
\end{split}
\end{equation}
Letting $a\rightarrow \infty,$  we obtain
\begin{align}\label{3.16}
 &\sum^{M}_{k=1}\bigg(\bigg.a_{k}-\sum^{l_{k}}_{l=1}b_{kl}e^{-\beta(1-\underline{\theta})t_{0}}\alpha_{kl}
 -\sum^{l_{k}}_{l=1}b_{kl}\frac{1}{\underline{\theta}}e^{-\beta(1-\underline{\theta})t_{0}}(1-\alpha_{kl})\bigg)\bigg.\mE\int^{t}_{t_0}U_{k}(x(s), s)\dif s\no\\
 &\leq c_{0} +a_{0}t.
 \end{align}
This means that  assertion (ii) holds.

(iii) Since $a_0=0,$ we derive from \eqref{f2} that
\begin{align*}
\mE[e^{\varepsilon t}U_{0}(x(t), t\wedge a, r(t))]
 \leq \bar{c}_{0}.
 \end{align*}
 This implies that \eqref{3c} holds.

 Using the similar method in (i) without taking the expectation, we can show that
\begin{align}
e^{\varepsilon t}U_{0}(x(t), t)\leq \bar{c}_{0} + M(t),
\end{align}
where $M(t)= \int^{t}_{0}e^{\varepsilon t}V_{x}(x(s),s, r(s))g(x_{s},t, r(s))\dif B(s).$
Due to Lemma 3.2,  it follows that
$$\limsup_{t\rightarrow \infty}e^{\varepsilon t}U_{0}(x(t), t)< \infty,\,a.s.$$
Thus, there exists a finite positive random variable $\eta$ such that
$$ \sup_{t_{0}\leq t <\infty}e^{\varepsilon t}U_{0}(x(t), t)<\eta,\,a.s.$$
Thus, the proof of  (\ref{3c1}) is complete.
\end{proof}

 We now illustrate the theoretical results in Theorem \ref{T3.3} by the following  example.
\begin{exa} {\rm
Let   $\nu(\cdot)$ be a probability measure on $[\underline{\theta}, 1]$.    Set $ S=\{1,2\}, \beta=0.5,  \lambda(\theta, t)=0.5(1-\theta), \underline{\theta}=0.5$, $d=1$. Let $r(t)$ be a Markov chain with generator
$$
\Gamma=\left(
  \begin{array}{ccc}
    -1 & 1 \\
    2 & -2 \\
  \end{array}
\right).
$$
 Consider the following equation:
\begin{align}
 &\dif x(t)= f( x_t, t,r(t))\dif t +g(  x_t, t, r(t))\dif  B(t), t\in [t_{0}, \infty)\no\\
 & x(t)=\xi(t), t\in [\underline{\theta} t_{0}, t_{0}],
 \end{align}
where for $\varphi \in \C$
\begin{align*}
f( \varphi, t,i)=
\begin{cases}
-5(\varphi(1)+\varphi^{3}(1)+\varphi^{5}(1)) +0.5\int^{1}_{\frac{1}{2}}e^{-0.5(1-\theta) t}|\varphi(\theta )|\dif \nu (\theta), i=1,\\
0.05\varphi(1) +0.05\int^{1}_{\frac{1}{2}}e^{-0.5(1-\theta) t}|\varphi(\theta )|\dif \nu (\theta), i=2;
\end{cases}
 \end{align*}
and
\begin{align*}
g(\varphi, t,i)
=\begin{cases}
0.5\int^{1}_{\frac{1}{2}}e^{-0.5(1-\theta) t}|\varphi(1)|^{2}|\varphi(\theta )|\dif \nu (\theta), i=1,\\
0.2\int^{1}_{\frac{1}{2}}e^{-0.5(1-\theta) t}|\varphi(\theta )|\dif \nu (\theta), i=2.
\end{cases}
 \end{align*}
Define
\begin{align*}
V( x, t,i)
=\begin{cases}
x^{2}, \quad i=1,\\
2(x^{2}+x^{6}), \quad  i=2.
\end{cases}
 \end{align*}
When $i=1, $  it follows that
\begin{align*}
&LV(\varphi, t,1)=2\varphi(1) f(\varphi, t,1)+|g(\varphi, t,1)|^{2}+ \sum^{2}_{j=1}\gamma_{1j}V(\varphi(1), t,j)\\
&\leq 2\varphi(1)[-5(\varphi(1)+\varphi^{3}(1)+\varphi^{5}(1))] +\varphi(1)\int^{1}_{\frac{1}{2}}e^{-0.5(1-\theta) t}|\varphi(\theta )|\dif \nu (\theta)\\
&\quad \quad \quad+0.25\int^{1}_{\frac{1}{2}}e^{-0.5(1-\theta) t}|\varphi(1)|^{4}|\varphi(\theta )|^{2}\dif \nu (\theta)-|\varphi(1)|^{2}+2(|\varphi(1)|^{2}+|\varphi(1)|^{6})\\
&\leq -9|\varphi(1)|^{2}-10|\varphi(1)|^{4}-8|\varphi(1)|^{6}+\varphi(1)\int^{1}_{\frac{1}{2}}e^{-0.5(1-\theta) t}|\varphi(\theta )|\dif \nu (\theta)\\
&\quad \quad \quad+0.25\int^{1}_{\frac{1}{2}}e^{-0.5(1-\theta) t}|\varphi(1)|^{4}|\varphi(\theta )|^{2}\dif \nu (\theta).
 \end{align*}
When $i=2, $  we have
\begin{align*}
&LV(\varphi, t,2)=(4\varphi(1)+12\varphi^{5}(1))f(\varphi, t,2)+0.02(4+60\varphi^{4}(1))\int^{1}_{\frac{1}{2}}e^{-0.5(1-\theta) t}|\varphi(\theta )|^{2}\dif \nu (\theta)\\
 &+ \sum^{2}_{j=1}\gamma_{2j}V(\varphi(1), t,j)\\
& \leq 0.2|\varphi(1)|^{2}+0.6|\varphi(1)|^{6}+0.2\varphi(1)\int^{1}_{\frac{1}{2}}e^{-0.5(1-\theta) t}|\varphi(\theta )|\dif \nu (\theta)\\
&+0.6\varphi^{5}(1)\int^{1}_{\frac{1}{2}}e^{-0.5(1-\theta) t}|\varphi(\theta )|\dif \nu (\theta)+0.08\int^{1}_{\frac{1}{2}}e^{-0.5(1-\theta) t}|\varphi(\theta )|^{2}\dif \nu (\theta)\\
&+1.2\int^{1}_{\frac{1}{2}}e^{-0.5(1-\theta) t}|\varphi(t)|^{4}|\varphi(\theta )|^{2}\dif \nu (\theta)+2|\varphi(1)|^{2}-4(|\varphi(1)|^{2}+|\varphi(1)|^{6})\\
& \leq -1.8|\varphi(1)|^{2}-3.4|\varphi(1)|^{6}\\
&+0.2\varphi(1)\int^{1}_{\frac{1}{2}}e^{-0.5(1-\theta) t}|\varphi(\theta )|\dif \nu (\theta)+0.6\varphi^{5}(1)\int^{1}_{\frac{1}{2}}e^{-0.5(1-\theta) t}|\varphi(\theta )|\dif \nu (\theta)\\
&+0.08\int^{1}_{\frac{1}{2}}e^{-0.5(1-\theta) t}|\varphi(\theta )|^{2}\dif \nu (\theta)+1.2\int^{1}_{\frac{1}{2}}e^{-0.5(1-\theta) t}|\varphi(1)|^{4}|\varphi(\theta )|^{2}\dif \nu (\theta)\\
& \leq -1.8|\varphi(1)|^{2}-3.4|\varphi(1)|^{6}\\
&+\int^{1}_{\frac{1}{2}}e^{-0.5(1-\theta) t}0.2(|\varphi(1)|^{2})^{\frac{1}{2}}(|\varphi(\theta )|^{2})^{\frac{1}{2}}\dif \nu (\theta)+\int^{1}_{\frac{1}{2}}e^{-0.5(1-\theta) t}0.6(|\varphi(1)|^{6})^{\frac{5}{6}}(|\varphi(\theta )|^{6})^{\frac{1}{6}}\dif \nu (\theta)\\
&+\int^{1}_{\frac{1}{2}}e^{-0.5(1-\theta) t}0.08|\varphi(\theta )|^{2}\dif \nu (\theta)+\int^{1}_{\frac{1}{2}}e^{-0.5(1-\theta) t}1.2(|\varphi(1)|^{6})^{\frac{4}{6}}(|\varphi(\theta )|^{6})^{\frac{2}{6}}\dif \nu (\theta)\\
& \leq -1.8|\varphi(1)|^{2}-3.4|\varphi(1)|^{6}\\
&+\int^{1}_{\frac{1}{2}}e^{-0.5(1-\theta) t}[0.2(|\varphi(1)|^{2})^{\frac{1}{2}}(|\varphi(\theta )|^{2})^{\frac{1}{2}}+0.08|\varphi(\theta )|^{2}]\dif \nu (\theta)\\
&+\int^{1}_{\frac{1}{2}}e^{-0.5(1-\theta) t}[0.6(|\varphi(1)|^{6})^{\frac{5}{6}}(|\varphi(\theta )|^{6})^{\frac{1}{6}}+1.2(|\varphi(1)|^{6})^{\frac{4}{6}}(|\varphi(\theta )|^{6})^{\frac{2}{6}}]\dif \nu (\theta).
\end{align*}
Then, we have
\begin{align*}
LV(\varphi, t,i) &\leq -1.8|\varphi(1)|^{2}-3.4|\varphi(1)|^{6}\\
&+\int^{1}_{\frac{1}{2}}e^{-0.5(1-\theta) t}[(|\varphi(1)|^{2})^{\frac{1}{2}}(|\varphi(\theta )|^{2})^{\frac{1}{2}}+0.08|\varphi(\theta )|^{2}]\dif \nu (\theta)\\
&+\int^{1}_{\frac{1}{2}}e^{-0.5(1-\theta) t}[0.6(|\varphi(1)|^{6})^{\frac{5}{6}}(|\varphi(\theta )|^{6})^{\frac{1}{6}}+1.2(|\varphi(1)|^{6})^{\frac{4}{6}}(|\varphi(\theta )|^{6})^{\frac{2}{6}}]\dif \nu (\theta).
 \end{align*}
Obviously,  we can choose
\begin{align*}
&U_{0}(x,t)=|x|^{2},U_{1}(x,t)=|x|^{2}, U_{2}(x,t)=|x|^{6}, a_{0}=0, a_{1}=1.8,\\
& a_{2}=3.4,, b_{11}=1, b_{12}=0.08,
b_{21}=0.6, b_{22}=1.2.
 \end{align*}
From Theorem 3.3,  we could know that the following results hold.
\begin{enumerate}
\item[(i)] $$\limsup_{t\rightarrow \infty}\mE[|x(t_{0},t,\xi,i_{0})|^{2}]=0,$$
\item[(ii)] $$\limsup_{t\rightarrow \infty}\frac{1}{t}\int^{t}_{t_{0}}\mE[|x(t_{0},s,\xi,i_{0})|^{2}]\dif s =0,$$
$$\limsup_{t\rightarrow \infty}\frac{1}{t}\int^{t}_{t_{0}}\mE[|x(t_{0},s,\xi,i_{0})|^{6}]\dif s =0,$$

\item[(iii)] \begin{align*}
\limsup_{t\rightarrow \infty}\frac{1}{t}\log (\mE[|x(t_{0},t,\xi,i_{0})|^{2}]) \leq -0.05.\no
 \end{align*}
  \begin{align*}
\limsup_{t\rightarrow \infty}\frac{1}{t}\log (|x(t_{0},t,\xi,i_{0})|^{2}) \leq -0.05,\,a.s.
 \end{align*}
\end{enumerate}
}
\end{exa}

 Now, we give the second example.

\begin{exa} {\rm
  Let $\nu_{1}(\cdot)$  be a  probability measure on $[\underline{\theta}, 1]$ and  $\nu_{2}(\cdot)=\delta_{1}(\cdot).$  Set $S=\{1,2\}, \beta=0.6,  \lambda(\theta, t)=0.6(1-\theta), \underline{\theta}=0.7, d=1$. Let $r(t)$ be a Markov chain with generator
$$
\Gamma=\left(
  \begin{array}{ccc}
    -1 & 1 \\
    3 & -3 \\
  \end{array}
\right).
$$
 Consider the following equation:
\begin{align}
 &\dif x(t)= f( x_t, t,r(t))\dif t +g(  x_t, t, r(t))\dif  B(t), t\in [t_{0}, \infty)\no\\
 & x(t)=\xi(t), t\in [\underline{\theta} t_{0}, t_{0}],
 \end{align}
where for $\varphi\in \C$
\begin{align*}
f( \varphi, t,i)=
\begin{cases}
-6(\varphi(1)+\varphi^{3}(1)+\varphi^{7}(1)) +\int^{1}_{0.7}e^{-0.6(1-\theta) t}\varphi(\theta )\dif \nu_{1} (\theta), i=1,\\
0.04\varphi(1) +0.04\int^{1}_{0.7}e^{-0.6(1-\theta) t}\varphi(\theta )\dif \nu_{2} (\theta), i=2;
\end{cases}
 \end{align*}
and
\begin{align*}
g( \varphi, t,i)
=\begin{cases}
0.5\int^{1}_{0.7}e^{-0.6(1-\theta) t}|\varphi(1)|^{2}|\varphi(\theta )|^{2}\dif \nu_{1} (\theta), i=1,\\
0.1\int^{1}_{0.7}e^{-0.6(1-\theta) t}\varphi(\theta )\dif \nu_{2} (\theta), i=2.
\end{cases}
 \end{align*}
From above equation, when $r(t)=2, $   by the definition of $\nu_2$ the equation becomes
\begin{align}
 &\dif x(t)= 0.08x(t)\dif t +0.1x(t)\dif  B(t), t\in [t_{0}, \infty)\no\\
 & x(t)=\xi(t), t\in [\underline{\theta} t_{0}, t_{0}],
 \end{align}
Obviously, the solution of the above equation will  blow up. But in the following, we will show that the overall system  is stable.
Set

\begin{align*}
V( x, t,i)
=\begin{cases}
x^{2}, i=1,\\
2x^{2}+3x^{8}, i=2.
\end{cases}
 \end{align*}
When $i=1, $  it follows that
\begin{align*}
&LV(\varphi, t,1)=2\varphi(1)f(\varphi, t,1)+|g(\varphi, t,1)|^{2}+ \sum^{2}_{j=1}\gamma_{1j}V(\varphi(1), t,j)\\
&\leq 2\varphi(1)[-6(\varphi(1)+\varphi^{3}(1)+\varphi^{7}(1))] +2|\varphi(1)|\int^{1}_{0.7}e^{-0.6(1-\theta) t}|\varphi(\theta )|\dif \nu_{1} (\theta)\\
&\quad \quad \quad+0.25\int^{1}_{0.7}e^{-0.6(1-\theta) t}|\varphi(1)|^{4}|\varphi(\theta )|^{4}\dif \nu_{1} (\theta)-|\varphi(1)|^{2}+2|\varphi(1)|^{2}+3|\varphi(1)|^{8}\\
&\leq -11|\varphi(1)|^{2}-12|\varphi(1)|^{4}-9|\varphi(1)|^{8}+2|\varphi(1)|\int^{1}_{0.7}e^{-0.6(1-\theta) t}|\varphi(\theta )|\dif \nu_{1} (\theta)\\
&\quad \quad \quad+0.25\int^{1}_{0.7}e^{-0.6(1-\theta) t}|\varphi(1)|^{4}|\varphi(\theta )|^{4}\dif \nu_{1} (\theta).
 \end{align*}
When $i=2, $  we have
\begin{align*}
&LV(\varphi, t,2)=(4\varphi(1)+24\varphi^{7}(1))f(\varphi, t,2)+0.005(4+168\varphi^{6}(1))\int^{1}_{0.7}e^{-0.6(1-\theta) t}|\varphi(\theta )|^{2}\dif \nu_{2} (\theta)\\
 &+ \sum^{2}_{j=1}\gamma_{2j}V(\varphi(1), t,j)\\
& \leq 0.08(4\varphi(1)+24\varphi^{7}(1))\varphi(1)+0.01(2+84\varphi^{6}(1))|\varphi(1)|^{2}+3\varphi^{2}(1)-3(2\varphi^{2}(1)+3\varphi^{8}(1))\\
& \leq -2.64|\varphi(1)|^{2}-6.24|\varphi(1)|^{8}.
\end{align*}
Then, we have
\begin{align*}
LV(\varphi, t,i) &\leq -2.66|\varphi(1)|^{2}-6.24|\varphi(1)|^{8}\\
&+2\varphi(1)\int^{1}_{0.7}e^{-0.6(1-\theta) t}|\varphi(\theta )|\dif \nu_{1} (\theta)
+0.25\int^{1}_{0.7}e^{-0.6(1-\theta) t}|\varphi(1)|^{4}|\varphi(\theta )|^{4}\dif \nu_{1} (\theta)\\
 &\leq -2.64|\varphi(1)|^{2}-6.24|\varphi(1)|^{8}\\
&+\int^{1}_{0.7}e^{-0.6(1-\theta) t}2(|\varphi(1)|^{2})^{\frac{1}{2}}(|\varphi(\theta )|^{2})^{\frac{1}{2}}\dif \nu_{1} (\theta)\\
&+\int^{1}_{0.7}e^{-0.6(1-\theta) t}0.25(|\varphi(1)|^{8})^{\frac{1}{2}}(|\varphi(\theta )|^{8})^{\frac{1}{2}}\dif \nu_{1} (\theta).
 \end{align*}
Obviously,  we can choose
\begin{align*}
&U_{0}(x,t)=|x|^{2},U_{1}(x,t)=|x|^{2}, U_{2}(x,t)=|x|^{8}, a_{0}=0, a_{1}=2.64,\\
& a_{2}=6.24,, b_{11}=2,
b_{21}=0.25.
 \end{align*}
From Theorem 3.3,  we have  the following results:
\begin{itemize}
\item[(i)] $$\limsup_{t\rightarrow \infty}\mE[|x(t_{0},t,\xi,i_{0})|^{2}]=0,$$
\item[(ii)] $$\limsup_{t\rightarrow \infty}\frac{1}{t}\int^{t}_{t_{0}}\mE[|x(t_{0},s,\xi,i_{0})|^{2}]\dif s =0,$$
$$\limsup_{t\rightarrow \infty}\frac{1}{t}\int^{t}_{t_{0}}\mE[|x(t_{0},s,\xi,i_{0})|^{8}]\dif s =0,$$

\item[(iii)] \begin{align*}
\limsup_{t\rightarrow \infty}\frac{1}{t}\log (\mE[|x(t_{0},t,\xi,i_{0})|^{2}]) \leq -0.1.\no
 \end{align*}
\begin{align*}
\limsup_{t\rightarrow \infty}\frac{1}{t}\log (|x(t_{0},t,\xi,i_{0})|^{2}) \leq -0.1,\,a.s.
 \end{align*}
\end{itemize}
}
\end{exa}

\subsection{Polynomial Stability}

In this subsection,  we will investigate  the polynomial stability of the solution for HPSFEDs \eqref{eq1}.

\bt
Assume (H1), and let (H2') hold with $a_0=0$. If
\begin{align}
-a_{k}+\sum^{l_{k}}_{l=1}b_{kl}\alpha_{kl}
 +\sum^{l_{k}}_{l=1}b_{kl}\frac{1}{\underline{\theta}}(1-\alpha_{kl})<0, \, k=1, 2, \cdots, M,
 \end{align}
 then
 the global solution $x(t_{0},t,\xi,i_{0})$   has almost surely polynomial stability, i.e.
\begin{align}\label{24}
\limsup_{t\rightarrow \infty }\frac{\log U_{0}(x(t),t)}{\log(1+t)}\leq -\varepsilon,
 \end{align}
where $\varepsilon$ is a positive constant satisfying $$-a_{k}+\sum^{l_{k}}_{l=1}b_{kl}\alpha_{kl}
 +\sum^{l_{k}}_{l=1}b_{kl}\underline{\theta}^{-(1+\varepsilon)}(1-\alpha_{kl})<0, k=2,3,\cdots, M,$$
and
$$
\varepsilon -a_{1}+\sum^{l_{1}}_{l=1}b_{1l}\alpha_{1l}
 +\sum^{l_{1}}_{l=1}b_{1l}\underline{\theta}^{-(1+\varepsilon)}(1-\alpha_{1l})<0.
$$
\et

\begin{proof}
 Define
$\sigma_{a}=\inf\{t\geq t_0: |x(t)|\geq a\}$
 as before.  Set
  \begin{align*}
M(t)=\int^{t}_{t_{0}}(1+s)^{\varepsilon}V_{x}(x(s), s, r(s))g(x_{s}, s,r(s))\dif B(s).
 \end{align*}
 Using It\^{o}'s formula and taking the expectation, we have
\begin{align}
&(1+t\wedge \sigma_{a})^{\varepsilon}V(x(t\wedge \sigma_a), t\wedge \sigma_a, r(t\wedge \sigma_a))]\no\\
&=(1+t_{0})^{\varepsilon}V(x(t_0), t_0, r(t_0))+ \int^{t\wedge  \sigma_a}_{t_{0}}\varepsilon (1+s)^{\varepsilon-1}V(x_{s},s,r(s))\dif s\no\\
& \quad \quad\quad+ \int^{t\wedge  \sigma_a}_{t_{0}}(1+s)^{\varepsilon}LV(x_{s},s,r(s))\dif s +M(t\wedge \sigma_{a}) \no\\
& \leq (1+t_{0})^{\varepsilon}V(x(t_0), t_0, r(t_0)) +\int^{t\wedge  \sigma_a}_{t_{0}}\varepsilon (1+s)^{\varepsilon}V(x_{s},s,r(s))\dif s \no\\
& \quad \quad\quad+  \int^{t\wedge  \sigma_a}_{t_{0}}(1+s)^{\varepsilon}\bigg\{\bigg. a_{0}+\sum^{M}_{k=1}\bigg[\bigg.-a_{k}U_{k}(x(s), s)\no\\
&\quad \quad \quad+ \sum^{l_{k}}_{l=1}b_{kl}\int^{1}_{\underline{\theta}} U_{k}^{\alpha_{kl}}(x(s), s) U_{k}^{1-\alpha_{kl}}(x(\theta s), \theta s)\dif \nu_{k}(\theta) \bigg]\bigg.\bigg\}\bigg.\dif s
+M(t\wedge \sigma_{a})\no\\
&\leq (1+t_{0})^{\varepsilon}V(x(t_0), t_0, r(t_0))+\int^{t\wedge  \sigma_a}_{t_{0}}\varepsilon (1+s)^{\varepsilon}V(x_{s},s,r(s))\dif s \no\\
& \quad \quad\quad+  \int^{t\wedge  \sigma_a}_{t_{0}}(1+s)^{\varepsilon}\bigg\{\bigg. a_{0}+\sum^{M}_{k=1}\bigg[\bigg.-a_{k}U_{k}(x(s), s)\no\\
&  \quad \quad\quad  +\sum^{l_{k}}_{l=1} b_{kl}\int^{1}_{\underline{\theta}} U_{k}^{\alpha_{kl}}(x(s), s) U_{k}^{1-\alpha_{kl}}(x(\theta s), \theta s)\dif \nu_{k}(\theta) \bigg]\bigg.\bigg\}\bigg.\dif s+M(t\wedge \sigma_{a})\no\\
&\leq (1+t_{0})^{\varepsilon}V(x(t_0), t_0, r(t_0)) +\int^{t\wedge  \sigma_a}_{t_{0}}\varepsilon (1+s)^{\varepsilon}V(x_{s},s,r(s))\dif s\no\\
 & \quad \quad\quad +  \int^{t\wedge  \sigma_a}_{t_{0}}(1+s)^{\varepsilon}\bigg\{\bigg. a_{0}+\sum^{M}_{k=1}\bigg[\bigg.-a_{k}U_{k}(x(s), s)+ \sum^{l_{k}}_{l=1}b_{kl}\alpha_{kl}U_{k}(x(s), s)\no\\
&  \quad \quad\quad +\sum^{l_{k}}_{l=1}b_{kl}(1-\alpha_{kl})\int^{1}_{\underline{\theta}} U_{k}(x(\theta s), \theta s)\dif \nu_{k}(\theta) \bigg]\bigg.\bigg\}\bigg.\dif s+M(t\wedge \sigma_{a})\no\\
&\leq \tilde{c}_{0} +\int^{t\wedge  \sigma_a}_{t_{0}}\varepsilon (1+s)^{\varepsilon}V(x_{s},s,r(s))\dif s\no\\
 & \quad \quad\quad +  \int^{t\wedge  \sigma_a}_{t_{0}}(1+s)^{\varepsilon}\bigg\{\bigg. \sum^{M}_{k=1}\bigg[\bigg.-a_{k}U_{k}(x(s), s)+ \sum^{l_{k}}_{l=1}b_{kl}\alpha_{kl}U_{k}(x(s), s)\no\\
&\quad\quad\quad+\sum^{l_{k}}_{l=1}b_{kl}\underline{\theta}^{-(1+\varepsilon)}(1-\alpha_{kl})U_{k}(\varphi( s),  s) \bigg]\bigg.\bigg\}\bigg.\dif s+M(t\wedge \sigma_{a}) \no\\
&\leq \tilde{c}_{0} +\int^{t_{0}}_{\underline{\theta}t_{0}}\sum^{M}_{k=1}\sum^{l_{k}}_{l=1}b_{kl}
\underline{\theta}^{-(1+\varepsilon)}(1-\alpha_{kl})U_{k}(\xi( s),  s)\dif s\no\\
 &+ \int^{t\wedge  \sigma_a}_{t_{0}}(1+s)^{\varepsilon}\bigg(\bigg.\varepsilon -a_{1}+\sum^{l_{1}}_{l=1}b_{1l}\alpha_{1l}
 +\sum^{l_{1}}_{l=1}b_{1l}\underline{\theta}^{-(1+\varepsilon)}(1-\alpha_{1l})\bigg)\bigg.U_{1}(x(s), s)\no\\
 &+\int^{t\wedge  \sigma_a}_{t_{0}}(1+s)^{\varepsilon}\bigg\{\bigg. \sum^{M}_{k=2}\bigg(\bigg.-a_{k}+\sum^{l_{k}}_{l=1}b_{kl}\alpha_{kl}
 +\sum^{l_{k}}_{l=1}b_{kl}\underline{\theta}^{-(1+\varepsilon)}(1-\alpha_{kl})\bigg)\bigg.U_{k}(x(s), s)\bigg\}\bigg.\dif s\no\\
 &\leq \tilde{c}_{0}+M(t\wedge \sigma_{a}),
 \end{align}
where $\tilde{c}_{0}=(1+t_{0})^{\varepsilon}V(x(t_0), t_0, r(t_0))+\int^{t_{0}}_{\underline{\theta}t_{0}}\sum^{M}_{k=1}\sum^{l_{k}}_{l=1}b_{kl}
\underline{\theta}^{-(1+\varepsilon)}(1-\alpha_{kl})U_{k}(\xi( s),  s)\dif s.$\\
By virtue of the conditions in the theorem, we have
\begin{align*}
(1+t)^{\varepsilon}U_{0}(x(t), t)
 \leq \tilde{c}_{0}+M(t).
 \end{align*}
Then,
\begin{align*}
\limsup_{t\rightarrow \infty}(1+t)^{\varepsilon}U_{0}(x(t), t)
 < \infty.
 \end{align*}
This implies the required assertion (\ref{24}) immediately.
 \end{proof}
 The following example illustrates the theory of polynomial stability.

\begin{exa} {\rm
  Let  $\nu_{1}(\cdot)$  be a  probability measure on $[\underline{\theta}, 1]$  and  $\nu_{2}(\cdot)=\delta_{1}(\cdot).$  Set $S=\{1,2\}, \underline{\theta}=0.75, d=1$. Let $r(t)$ be a Markov chain with generator
$$
\Gamma=\left(
  \begin{array}{ccc}
    -1 & 1 \\
    4 & -4 \\
  \end{array}
\right).
$$
 Consider the following HPSFDE:
\begin{align}
 &\dif x(t)= f( x_t, t,r(t))\dif t +g(  x_t, t, r(t))\dif  B(t), t\in [t_{0}, \infty)\no\\
 & x(t)=\xi(t), t\in [\underline{\theta} t_{0}, t_{0}],
 \end{align}
where for $\varphi \in \C$
\begin{align*}
f( \varphi, t,i)=
\begin{cases}
-6(\varphi(1)+\varphi^{3}(1)+\varphi^{7}(1)) +0.5\int^{1}_{0.75}\varphi(\theta )\dif \nu_{1} (\theta), r(t)=1,\\
0.04\varphi(1) +0.03\int^{1}_{0.75}\varphi(\theta )\dif \nu_{2} (\theta), r(t)=2;
\end{cases}
 \end{align*}
and
\begin{align*}
g( \varphi, t,i)
=\begin{cases}
0.2\int^{1}_{0.75}|\varphi(1)|^{1.5}|\varphi(\theta )|^{2.5}\dif \nu_{1} (\theta), r(t)=1,\\
0.1\int^{1}_{0.75}|\varphi(\theta )|\dif \nu_{2} (\theta), r(t)=2.
\end{cases}
 \end{align*}
From above equation, when $r(t)=2, $ we can see that  the equation is
\begin{align}
 &\dif x(t)= 0.07x(t)\dif t +0.1x(t)\dif  B(t), t\in [t_{0}, \infty)\no\\
 & x(t)=\xi(t), t\in [\underline{\theta} t_{0}, t_{0}],
 \end{align}
obviously, the solution of the above equation will  blow up. But in the following, we will show that the overall system  is polynomial stable.
Set
\begin{align*}
V( x, t,i)
=\begin{cases}
x^{4}, i=1,\\
2x^{4}+3x^{10}, i=2.
\end{cases}
 \end{align*}
When $i=1, $  it follows that
\begin{align*}
&LV(\varphi, t,1)=4\varphi^{3}(1)f(\varphi, t,1)+0.24\varphi^{2}(1)|g(\varphi, t,1)|^{2}+ \sum^{2}_{j=1}\gamma_{1j}V(\varphi(1), t,j)\\
&\leq 4\varphi^{3}(1)[-6(\varphi(1)+\varphi^{3}(1)+\varphi^{7}(1))] +2\varphi^{3}(1)\int^{1}_{0.75}|\varphi(\theta )|\dif \nu_{1} (\theta)\\
&\quad \quad \quad \quad\quad \quad \quad+0.24\int^{1}_{0.75}|\varphi(1)|^{5}|\varphi(\theta )|^{5}\dif \nu_{1} (\theta)-|\varphi(1)|^{4}+2|\varphi(1)|^{4}+3|\varphi(1)|^{10}\\
&\leq -23|\varphi(1)|^{4}-24|\varphi(1)|^{6}-21|\varphi(1)|^{10}+2|\varphi(1)|^{3}\int^{1}_{0.75}|\varphi(\theta )|\dif \nu_{1} (\theta)\\
&\quad \quad \quad \quad\quad \quad \quad+0.24\int^{1}_{0.75}|x(t)|^{5}|\varphi(\theta )|^{5}\dif \nu_{1} (\theta).
 \end{align*}
When $i=2, $  we have
\begin{align*}
&LV(\varphi, t,2)=(8\varphi^{3}(1)+30\varphi^{9}(1))f(\varphi, t,2)+0.005(24\varphi^{2}(1)+270\varphi^{8}(1))\int^{1}_{0.75}|\varphi(\theta )|^{2}\dif \nu_{2} (\theta)\\
 &+ \sum^{2}_{j=1}\gamma_{2j}V(x(t), t,j)\\
& \leq (8\varphi^{3}(1)+30\varphi^{9}(1))0.07\varphi(1)+0.01(12\varphi^{2}(1)+135\varphi^{8}(1))\varphi^{2}(1)\\
&+4\varphi^{4}(1)-4(2\varphi^{4}(1)+3\varphi^{10}(1))\\
& \leq -3.32|\varphi(1)|^{4}-8.55|\varphi(1)|^{10}.
\end{align*}
Then, we have
\begin{align*}
LV(\varphi, t,i) &\leq -3.32|\varphi(1)|^{4}-8.55|\varphi(1)|^{10}\\
&+2|\varphi(1)|^{3}\int^{1}_{0.75}|\varphi(\theta )|\dif \nu_{1} (\theta)
+0.24\int^{1}_{0.75}|x(t)|^{5}|\varphi(\theta )|^{5}\dif \nu_{1} (\theta)\\
 &\leq -3.32|\varphi(1)|^{2}-8.25|\varphi(1)|^{8}\\
&+\int^{1}_{0.75}2(|\varphi(1)|^{4})^{\frac{3}{4}}(|\varphi(\theta )|^{4})^{\frac{1}{4}}\dif \nu_{1} (\theta)\\
&+\int^{1}_{0.75}0.24(|\varphi(1)|^{10})^{\frac{1}{2}}(|\varphi(\theta )|^{10})^{\frac{1}{2}}\dif \nu_{1} (\theta)
 \end{align*}
Obviously,  we can choose
\begin{align*}
&U_{0}(x,t)=|x|^{4},U_{1}(x,t)=|x|^{4}, U_{2}(x,t)=|x|^{10}, a_{0}=0, a_{1}=3.32,\\
& a_{2}=8.55,, b_{11}=2,
b_{21}=0.24.
 \end{align*}
From Theorem 3.6,  we conclude that the overall system  is  polynomial stable.
}
Obviously, this example is similar to example 3.5, but $f,g$ in this example satisfy  (H2)' while $f,g$ in  example 3.5 satisfy  (H2). This  difference leads to different stable properties of the solution.
\end{exa}

\section*{Acknowledgements}
This research was supported by the National Natural Science Foundation of China (Grant no.61876192, 11626236), the Fundamental Research Funds for the Central Universities of South-Central University for Nationalities (Grant nos. CZY15017, KTZ20051, CZT20020)

\end{document}